\documentclass[12pt]{amsart}
\usepackage{latexsym}
\numberwithin{equation}{section}
\makeatletter
\renewcommand{\subsection}{\@startsection
{subsection}{2}{0mm}{\baselineskip}{-0.25cm}
{\normalfont\normalsize\rmfamily}}
\makeatother
\newtheorem{theorem}{Theorem}[section]
\newtheorem{lemma}[theorem]{Lemma}

\newtheorem{corollary}[theorem]{Corollary}
\newtheorem{theorem'}{Theorem}[subsection]
\newtheorem{lemma'}[theorem']{Lemma}
\theoremstyle{definition}
\newtheorem{definition}[theorem]{Definition}
\newtheorem{remark}[theorem]{Remark}
\newtheorem{remark'}[theorem']{Remark}
\newenvironment{claim*}{\begin{trivlist}\item[\hskip%
\labelsep{\itshape{Claim}.}]\it }%
{\end{trivlist}}
\def\[x]{[(g+1)/2]}
\begin{document}
\author[F.~Torres]{Fernando Torres}\thanks{The author is supported by a
grant from the International Atomic Energy Agency and UNESCO}
\title[Numerical semigroups]{Remarks on numerical semigroups}
\address{ICTP, Mathematics Section, P.O. Box 586, 34100, Trieste - Italy}
\email{feto@ictp.trieste.it}
\begin{abstract}
We extend results on Weierstrass semigroups at ramified points of double
covering of curves to any numerical semigroup whose genus is large
enough. As an application we strengthen the properties concerning
Weierstrass weights stated in \cite{To}.
\end{abstract}
\maketitle
\setcounter{section}{-1}
\section{Introduction}
Let $H$ be a numerical semigroup, that is, a subsemigroup of $(\mathbb N,
+)$ whose complement is finite. Examples of such
semigroups are  the Weierstrass semigroups at non-singular points of
algebraic curves.

In this paper we deal with the following type of semigroups:
\begin{definition}\label{def} Let $\gamma\ge 0$ an integer. $H$ is
called $\gamma$-hyperelliptic if the following conditions hold:
\begin{itemize}
\item[($E_1$)] $H$ has $\gamma$ even elements in $[2,4\gamma]$.
\item[($E_2$)] The $(\gamma+1)$th positive element of $H$ is $4\gamma+2$.
\end{itemize}
A 0-hyperelliptic semigroup is usually called hyperelliptic.
\end{definition}

The motivation for study of such semigroups comes from the study of
Weierstrass semigroups at ramified points of double coverings
of curves. Let $\pi: X\to \tilde X$ be a double covering of projective,
irreducible, non-singular algebraic curves over an algebraically closed
field $k$. Let $g$ and $\gamma$ be the genus of $X$ and
$\tilde X$ respectively. Assume that there exists $P\in X$ which is ramified
for $\pi$, and  denote by
$m_i$ the $i$th non-gap at $P$. Then the
Weierstrass semigroup $H(P)$ at $P$ satisfies the
following properties (cf. \cite{To}, \cite[Lemma 3.4]{To1}):
\begin{itemize}
\item[($P_1$)] $H(P)$ is $\gamma$-hyperelliptic, provided $g\ge
4\gamma+1$ if ${\rm char}(k)\neq 2$, and $g\ge 6\gamma-3$ otherwise.
\item[($P_2$)] $m_{2\gamma+1}=6\gamma+2$, provided $g\ge 5\gamma+1$.
\item[($P_3$)] $m_{\frac{g}{2}-\gamma-1}=g-2$ or
$m_{\frac{g-1}{2}-\gamma}=g-1$, provided $g\ge 4\gamma+2$.
\item[($P_4$)] The weight $w(P)$ of $H(P)$ satisfies
$$
\binom{g-2\gamma}{2}\le w(P)< \binom{g-2\gamma+2}{2}.
$$
\end{itemize}
Conversely if $g$ is large enough and if any of the above
properties is satisfied, then $X$ is a double covering of a curve of
genus $\gamma$. Aposteriori the
four above properties become equivalent whenever $g$ is large enough.

The goal of this paper is to extend these results for any
semigroup $H$ such that $g(H):= \#(\mathbb N\setminus H)$ is large enough.
We remark that there exist semigroups of genus large enough that cannot be
realized as
Weierstrass semigroups (see \cite{Buch1}, \cite[Scholium 3.5]{To}).

The key tool used to prove these equivalences is Theorem 1.10 in Freiman's
book \cite{Fre} which have to do with addition of finite sets. From this
theorem we deduce Corollary \ref{cor-cast} which can be considered as
analogous to Castelnuovo's genus bound for
curves in projective spaces (\cite{C}, \cite[p.116]{ACGH},
\cite[Corollary 2.8]{R}). Castelnuovo's result is
the key tool to deal with Weierstrass semigroups. This Corollary can also be
proved by means of properties of addition of residue classes (see
Remark \ref{cauchy}).

In \S2 we prove the equivalences $(P_1)\Leftrightarrow
(P_2)\Leftrightarrow (P_3)$. The equivalence $(P_1)\Leftrightarrow
(P_2)$ is proved under
the hypothesis $g(H)\ge 6\gamma+4$, while $(P_1)\Leftrightarrow (P_3)$
is proved under $g(H)=6\gamma+5$ or $g(H)\ge 6\gamma+8$. In both cases
the bounds on $g(H)$ are sharp (Remark \ref{sharp}). We mention that the
cases $\gamma\in\{1,2\}$ of $(P_1)\Leftrightarrow (P_3)$ were fixed by
Kato \cite[Lemmas 4,5,6,7]{K2}.

In \S3 we deal with the equivalence $(P_1)\Leftrightarrow (P_4)$. To
this purpose we determine bounds for the weight $w(H)$ of
the semigroup $H$, which is defined by
$$
w(H):= \sum_{i=1}^{g}(\ell_i -i),
$$
where $g:=g(H)$ and $\mathbb N\setminus H = \{\ell_1,\ldots,\ell_g\}$. It is
well
known that $0\le w(H)\le \binom{g}{2}$; clearly $w(H)=0\Leftrightarrow
H=\{g+i:i\in \mathbb Z^+\}$, and one can show that $w(H)=\binom{g}{2}
\Leftrightarrow H$ is $\mathbb N$, or $g(H)\ge 1$ and $H$ is hyperelliptic
(see e.g. \cite[Corollary III.5.7]{F-K}).
Associated to $H$ we have the number
\begin{equation}\label{even-gap}
\rho=\rho(H):=\{\ell \in \mathbb N\setminus H: \ell\ {\rm even}\}.
\end{equation}
In \cite[Lemma 2.3]{To} it has been shown that $\rho(H)$ is the unique
number $\gamma$ satisfying $(E_1)$ of Definition \ref{def}, and
\begin{itemize}
\item[($E_2'$)] $4\gamma+2 \in H$.
\end{itemize}
Thus we observe the following:
\begin{lemma}\label{feto0} Let $H$ be a $\gamma$-hyperelliptic
semigroup. Then
$$
\rho(H)=\gamma.
$$
\end{lemma}
We also observe that if $g(H)\ge 1$, then $H$ is hyperelliptic if and
only if $\rho(H)=0$. In general, $\rho(H)$ affects
the values of $w(H)$. Let us assume that $\rho(H)\ge
1$ (hence $w(H)<\binom{g}{2}$); then we find
$$
\binom{g-2\rho}{2}\le w(H)\le \left\{
\begin{array}{ll}
\binom{g-2\rho}{2}+2\rho^2 &  {\rm if\ } g\ge 6\rho+5 \\
\frac{g(g-1)}{3}   &  {\rm otherwise}
\end{array}
\right.
$$
(see Lemmas \ref{bo-weight} and \ref{opt-weight}). These bounds
strengthen results of Kato \cite[Thm.1]{K1} and Oliveira
\cite[p.435]{Oliv} (see Remark \ref{oliv}). From this
result we prove $(P_1)\Leftrightarrow (P_4)$ (Theorem
\ref{char-weight1}) under the hypothesis
\begin{equation}\label{bound-g}
g(H)\ge \left\{
\begin{array}{ll}
{\rm max}\{12\gamma-1,1\} & {\rm if\ } \gamma\in\{0,1,2\}, \\
11\gamma+1   & {\rm if\ } \gamma \in\{3,5\}, \\
\frac{21(\gamma-4)+88}{2}  & {\rm if\ }\gamma \in \{4,6\}, \\
\gamma^2+4\gamma+3   & {\rm if\ } \gamma\ge 7.
\end{array}
\right.
\end{equation}
The cases $\gamma \in \{1,2\}$ of that
equivalence was fixed by Garcia (see
\cite{G}). In this section we use ideas from Garcia's \cite[Proof of
Lemma 8]{G} and Kato's \cite[p. 144]{K1}.

In \S1 we recollet some arithmetical properties of numerical semigroups.
We mainly remark the influence of $\rho(H)$ on $H$.

It is a pleasure to thank Pablo Azcue and Gustavo T. de A. Moreira for
discussions about \S2.
\section{Preliminaries}
Throughout this paper we use the following notation
\begin{itemize}
\itemsep=0.5pt
\item semigroup:\quad numerical semigroup.
\item Let $H$ be a semigroup. The {\it genus} of $H$ is the number
$g(H):= \#(\mathbb N\setminus H)$, which throughout this article will be
assumed  bigger than 0. The positive elements of $H$ will be
called the {\it non-gaps} of $H$, and those of $G(H):= \mathbb N\setminus H$
will be called the {\it gaps} of $H$. We denote by $m_i(H)$ the $i$th
non-gap of $H$. If a semigroup is generated by $m,n,\ldots $ we denote
$H=\langle m,n,\ldots \rangle$.
\item $[x]$ stands for the integer part of $x\in \mathbb R$.
\end{itemize}
In this section we recall some arithmetical properties of semigroups. Let
$H$ be a semigroup of genus $g$. Set $m_j:= m_j(H)$ for each $j$. If
$m_1=2$ then $m_i=2i$ for $i=1,\ldots,g$. Let $m_i\ge 3$. By
the semigroup property of $H$ the first $g$ non-gaps satisfy the
following inequalities:
\begin{equation}\label{prop-sem}
m_i\ge 2i+1\ \ {\rm for}\ i=1,\ldots,g-2,\ \ m_{g-1}\ge 2g-2,\ \ m_g=2g
\end{equation}
(see \cite{Buch}, \cite[Thm.1.1]{Oliv}).
\medskip

Let $\rho$ be as in (\ref{even-gap}). From \cite[Lemma
2.3]{To1} we have that
\begin{equation}\label{feto}
\{4\rho+2i: i\in \mathbb N\} \subseteq H.
\end{equation}
{}From the definition of $\rho$, $H$ has
$\rho$ odd non-gaps in $[1,2g-1]$. Let denote
by
$$
u_{\rho} <\ldots < u_1
$$
such non-gaps.
\begin{lemma}\label{feto1} Let $H$ be a semigroup of genus $g$, and
$\rho$ the number of even gaps of $H$. Then
$$
2g \ge 3\rho.
$$
\end{lemma}
\begin{proof}
Let us assume that $g\le 2\rho -1$. From $u_1\le 2g-1$ we have
$u_{2\rho-g+1}\le 4g-4\rho-1$. Let $\ell$ be the biggest even
gap of $H$. Then $\ell \le 4g-4\rho$. For suppose that $\ell \ge
4g-4\rho+2$. Thus $\ell-u_{2\rho-g+1}\ge 3$,
and then $H$ would has $g-\rho+1$ odd gaps, namely $1,
\ell-u_{2\rho-g+1},\ldots, \ell-u_{\rho}$, a contradicition. Now
since in $[2,4g-4\rho]$ there are $2g-2\rho$ even numbers such that
$\rho$ of them are gaps, the lemma follows.
\end{proof}

Denote by $f_i:=f_i(H)$ the $i$th even non-gap of $H$.
Hence by (\ref{feto}) we have
\begin{equation}\label{gene}
H=\langle f_1,\ldots,
f_{\rho},4\rho+2,u_{\rho},\ldots,u_1\rangle.
\end{equation}
Observe that $f_{\rho}=4\rho$, and
\begin{equation}\label{even}
f_{g-\rho}=2g.
\end{equation}
By \cite[Lemma 2.1]{To1} and since $g\ge 1$ we have
\begin{equation}\label{des-odd-1}
u_{\rho} \ge\ {\rm max}\{2g-4\rho +1, 3\}.
\end{equation}
In particular, if $g\ge 4\rho$ we obtain
\begin{equation}\label{first=even}
m_1=f_1, \ldots, m_{\rho}=f_{\rho}.
\end{equation}
Note that (\ref{des-odd-1}) is only meanful for $g\ge 2\rho$. For
$g\le2\rho-1$ we have:
\begin{lemma}\label{feto2}
Let $H$ be a semigroup of genus $g$, and $\rho$ the number of even
gaps of $H$. If $g\le 2\rho-1$, then
$$
u_{\rho}\ge 4\rho -2g+1.
$$
\end{lemma}
\begin{proof}
{}From the proof of Lemma \ref{feto1} we have that $H$ has $2g-3\rho$ even
non-gaps in $[2,4g-4\rho]$. Consider the following sequence of even
non-gaps:
$$
2u_{\rho}<\ldots< u_{\rho}+u_{4\rho-2g}.
$$
Since in this sequence we have $2g-3\rho+1$ even non-gaps, then
$$
u_{\rho}+u_{4\rho-2g}\ge 4g-4\rho+2.
$$
Now, since $u_{4\rho-2g}\le 6g-8\rho+1$ the proof follows.
\end{proof}
%
%
%
\section{$\gamma$-hyperelliptic semigroups}
In this section we deal with properties $(P_1)$, $(P_2)$ and
$(P_3)$ stated in \S0. For $i\in \mathbb Z^+$ set
$$
d_i(H):= \gcd(m_1(H),\ldots,m_i(H)).
$$
\begin{theorem}\label{char1} Let $\gamma \in \mathbb N$, $H$
a semigroup of genus $g \ge 6\gamma +4$ if $\gamma\ge 1$. Then the
following statements are equivalent:
\begin{itemize}
\item[(i)] $H$ is $\gamma$-hyperelliptic.
\item[(ii)] $m_{2\gamma+1}(H)= 6\gamma +2$.
\end{itemize}
\end{theorem}
\begin{theorem}\label{char2}  Let $\gamma$ and $H$ be as in Theorem
\ref{char1}, and assume that $g\ge
1$ if $\gamma=0$. Then the following statements are equivalent:
\begin{itemize}
\item[(i)] $H$ is $t$-hyperelliptic for some $t\in \{0,\ldots,\gamma\}$.
\item[(ii)] $m_{2\gamma+1}(H)\le 6\gamma+2$.
\item[(iii)] $\rho(H)\le \gamma$.
\end{itemize}
\end{theorem}
\begin{theorem}\label{char3} Let $\gamma \in \mathbb N$, $H$ a semigroup of
genus $g=6\gamma+5$ or
$g\ge 6\gamma+7$. Set $r:= \[x]-\gamma-1$. Then the following statements
are equivalent:
\begin{itemize}
\item[(i)] $H$ is $\gamma$-hyperelliptic.
\item[(ii)] $m_r(H)=g-2$ if $g$ is even; $m_r(H)=g-1$ if $g$ is odd.
\item[(iii)] $m_r(H)\le g-1 < m_{r+1}(H)$.
\end{itemize}
\end{theorem}
\begin{theorem}\label{char4} Let $\gamma$, $H$ and $r$ be as in Theorem
\ref{char3}. Then the  following statements are equivalent:
\begin{itemize}
\item[(i)] $H$ is $t$-hyperelliptic for some $t\in\{0,\ldots,\gamma\}$.
\item[(ii)] $m_r(H)\le g-2$ if $g$ is even;
$m_r(H)\le g-1$ if $g$ is odd.
\item[(iii)] $m_r(H)\le g-1$.
\item[(iv)] $\rho(H)\le \gamma$.
\end{itemize}
\end{theorem}

To prove these results we need a particular case of the result below.
For $K$ a subset of a group we set $2K:=\{a+b: a, b \in K\}$.
\begin{lemma}[[Fre, Thm. 1.10 {]}]\label{thm-fre} Let
$K=\{0<m_1<\ldots<m_i\}
\subseteq \mathbb Z$ be such that $\gcd(m_1,\ldots,m_i)=1$. If $m_i\ge
i+1+b$, where $b$ is an integer satisfying $0\le b<i-1$, then
$$
\# 2K \ge 2i+2+b.
$$
\hfill $\Box$
\end{lemma}
\begin{corollary}\label{cor-cast} Let $H$ be a semigroup of genus
$g$, and $i\in \mathbb Z^+$. If
$$
d_i(H)= 1\qquad{\rm and}\qquad i\le g+1,
$$
then we have
$$
2m_i(H) \ge m_{3i-1}(H).
$$
\end{corollary}
\begin{proof} Let $K:=\{0,m_1(H),\ldots,m_i(H)\}$. Then by
(\ref{prop-sem}), we can apply Lemma
\ref{thm-fre} to $K$ with $b=i-2$.
\end{proof}
\begin{remark}\label{rem-cast} Both the hypothesis $d_i(H)=1$ and $i\le
g+1$ of the corollary above are necessaries. Moreover the conclusion of
that result is sharp:
\begin{itemize}
\item[(i)] Let $i=g+h$, $h\ge 2$. Then $2m_{g+h}=m_{3i-h}$.
\item[(ii)] Let $m_1=4$, $m_2=6$ and $m_3=8$. Then $d_3=2$ and $2m_3=m_7$.
\item[(iii)] Let $m_1=5$, $m_2=10$, $m_3=15$, $m_4=18$,
$m_5=20$. Then $2m_6=m_{14}$.
\end{itemize}
\end{remark}
\begin{remark}\label{cauchy} (i) The Corollary above can also be proved by
using
results on the addition of residue classes: let $H$ and $i$ be as in
\ref{cor-cast}; assume further that $2\le i\le g-2$ (the remaining cases
are easy to prove), and consider
$\tilde K:=\{m_1,\ldots,m_i\}\subseteq \mathbb Z_{m_i}$ (i.e. a subset of
the
integers modulus $m_i$). Let $N:= \# 2\tilde K$. Then it is easy
to see that
$$
2m_i\ge m_{i+N}.
$$
Consequently we have a proof of the above Corollary provided $N\ge
2i-1\ (*)$. Since $m_i\ge 2i+1$ (see (\ref{prop-sem})), we get
$(*)$ provided $m_i$ prime (Cauchy \cite{Dav1}, Davenport \cite{Dav},
\cite[Corollary 1.2.3]{M}), or provided
$\gcd(m_j,m_i)=1$ for $j=1,\ldots,i-1$ (Chowla \cite[Satz 114]{Lan},
\cite[Corollary 1.2.4]{M}). In
general by using the hypothesis $d_i(H)=1$ we can reduce the proof of the
Corollary to the case $\gcd(m_{i-1},m_i)=1$. Then we apply Pillai's
\cite[Thm 1]{Pi} generalization of Davenport and Chowla results (or
Mann's result \cite[Corollary 1.2.2]{M}).

(ii) Let $H$ and $i$ be as above and assume that
$2m_i=m_{3i-1}$. Then from (i) we have $N=\# 2\tilde K=2i-1$. Thus by
Kemperman \cite[Thm 2.1]{Kem} (or by \cite[Thm. 1.11]{Fre}) $2\tilde K$
satisfies one of the following conditions: (1) there exist $m, d\in \mathbb
Z_{m_i}$, such that $2\tilde K=\{m+dj:j=0,1,\ldots,N-1\}$, or (2) there
exists
a subgroup $F$ of $\mathbb Z_{m_i}$ of order $\ge 2$, such that $2\tilde K$
is the disjoint union of a non-empty set $C$ satisfying $C+F=C$, and a
set $C'$ satisfying $C'\subseteq c+F$ for some $c\in C'$. For instance
example (iii) of \ref{rem-cast} satisfies property (2).
\end{remark}

Set $m_j:= m_j(H)$ for each $j$.
\begin{proof} {\it (Theorem \ref{char1}).} By definition
$H$ is hyperelliptic if and only if $m_1=2$.
So let us assume that $\gamma\ge 1$.

(i) $\Rightarrow$ (ii): From Lemma \ref{feto0} and (\ref{des-odd-1}) we
find that $u_{\gamma}\ge 6\gamma+3$ if
$g\ge 5\gamma+1$. Then (ii) follows from (\ref{first=even}) and (\ref{feto}).

(ii) $\Rightarrow$ (i): We claim that $d_{2\gamma+1}(H)=2$. For
suppose that $d_{2\gamma+1}(H) \ge 3$. Then $6\gamma+2= m_{2\gamma+1}
\ge m_1+ 6\gamma$ and so $m_1\le 2$, a contradiction. Hence
$d_{2\gamma+1}(H)\le
2$. Now suppose that $d_{2\gamma+1}(H) =1$. Then Corollary \ref{cor-cast}
implies $$
2(6\gamma+2) = 2m_{2\gamma+1} \ge m_{6\gamma +2}.
$$
But, since $g-2 \ge 6\gamma +2$, by (\ref{prop-sem}) we
would have
$$
m_{6\gamma +2} \ge 2(6\gamma+2) +1
$$
which leads to a contradiction. This shows that $d_{2\gamma+1}(H)=2$. Now
since $m_{2\gamma+1}=6\gamma+2$ we have that $m_\gamma \le 4\gamma$.
Moreover,  there exist $\gamma$ even gaps of
$H$ in $[2, 6\gamma+2]$. Let $\ell$ be an even gap of $H$. The proof
follows from the following claim:
\begin{claim*}
$\ell <m_\gamma$.
\end{claim*}
\begin{proof} {\it (Claim).} Suppose that there exists an even gap
$\ell$ such that $\ell>m_\gamma$. Take the smallest $\ell$ with such a
property; then the following $\gamma$ even gaps: $\ell-m_\gamma<
\ldots,
\ell -m_1$ belong to $[2,m_\gamma]$. Thus, since $m_1>2$, we must have
$\ell-m_\gamma=2$. This implies
that $H$ has $\gamma+1$ even
non-gaps in $[2,6\gamma+2]$, namely $\ell-m_\gamma,\ldots,\ell-m_1,\ell$,
a contradiction.
\end{proof}
This finish the proof of Theorem \ref{char1}.
\end{proof}
\begin{proof} {\it (Theorem \ref{char2}).} The case $\gamma=0$
is trivial; so let assume $\gamma\ge 1$.

(i) $\Rightarrow$ (ii): Since $g\ge
5\gamma+1\ge 5t+1$ by Theorem \ref{char1} we have $m_{2t+1}=6t+2$.
Thus (ii) follows from Lemma \ref{feto0} and (\ref{feto}).

(ii) $\Rightarrow$ (iii): From the proof of (ii) $\Rightarrow$
(i) of Theorem \ref{char1} it follows that $d_{2\gamma+1}(H)=2$.
Consequently by using the hypothesis on $m_{2\gamma+1}$, and again from
the mentioned proof we have that all the gaps of $H$ belong to
$[2,m_\gamma]$.
Since $m_\gamma\le 4\gamma$ then we have $\rho(H)\le \gamma$

(iii) $\Rightarrow$ (i) Since $g\ge 4\gamma+1\ge 4\rho(H)+1$, the
proof follows from $(E_1)$ and $(E_2')$ (see \S0).
\end{proof}
\begin{proof} {\it (Theorem \ref{char3}).} (i) $\Rightarrow$
(ii): Similar to the proof of (i)
$\Rightarrow$ (ii) of Theorem \ref{char1} (here we need $g\ge
4\gamma+3$ (resp. $g\ge 4\gamma+4$) if $g$ is odd (resp. even)).
\smallskip

Before proving the other implications we remark that $g\le 3r-1$: in
fact, if $g\ge 3r$ we would have $g\le 6\gamma+6$ (resp. $g\ge
6\gamma+3$) provided $g$ even (resp. odd) - a contradiction.
\smallskip

(ii) $\Rightarrow$ (iii): Let $g$ even and suppose that $m_{r+1}=g-1$.
Then by
Corollary \ref{cor-cast} we would have $2g-2=2m_{r+1}\ge m_{3r+2}$ and
hence $g-1\ge 3r+2$. This contradicts the previous remark.

(iii) $\Rightarrow$ (i): We claim that $d_r(H)= 2$. Suppose that
$d_r(H)\ge
3$. Then we would have $g-1\ge m_r \ge m_1+3(r-1) \ge 3r-1$, which
contradicts the previous remark. Now suppose that
$d_r(H)=1$. Then by Corollary \ref{cor-cast} we would have
$$
2g-2\ge m_r \ge m_{3r-1},
$$
which again contradicts the previous remark.

Thus the number of even gaps of $H$ in $[2,g-1]$ is $\gamma$, and
$m_\gamma\le 4\gamma$. Let
$\ell$ be an even gap of $H$. As in the proof of the Claim in Theorem
\ref{char1} here we also have that $\ell<m_\gamma$. Now the proof
follows.
\end{proof}
\begin{proof} {\it (Theorem \ref{char4}).} (i) $\Rightarrow$
(ii): By Theorem \ref{char3} and since $t\le \gamma$ we have $g-2
=m_{g/2-t-1}$ or $g-1=m_{(g-1)/2-t}$. This implies (ii). The
implication (ii) $\Rightarrow$ (iii) is trivial.

(iii) $\Rightarrow$ (iv): As in the proof of Theorem \ref{char3} we
obtain $d_r(H)=2$. Then the number of even gaps of $H$ in $[2,g-1]$ is at
most $\gamma$. We claim that all the even gaps of
$H$ belong to that interval. For suppose there exists an even gap
$\ell>g-1$. Choose $\ell$
the smallest one and consider the even gaps $\ell-m_1<\ldots<\ell-m_r\le
g-1$. Then $r\le \gamma$ which yields to $g\le 4\gamma+2$,
a contradiction. Consequently $\rho(H)\le \gamma$.

The implication (iv) $\Rightarrow$ (i) follows from Theorem
\ref{char2}.
\end{proof}
\begin{remark}\label{sharp} The hypothesis on the genus in the
above theorems is sharp. To see this let $\gamma\ge 0$ an
integer, and let $X$ be the curve defined by the equation
$$
y^4=\mathop{\prod}\limits^{I}_{j=1}(x-a_j),
$$
where the $a_j's$ are pairwise distinct elements of a field $k$,
$I=4\gamma+3$ if
$\gamma$ is odd; $I=4\gamma+5$ otherwise. Let $P$ be the unique point over
$x=\infty$. Then $H(P)=\langle 4, I\rangle$ and so
$g(H(P))=6\gamma+3$ (resp. $6\gamma+6$), $m_{2\gamma+1}(H(P))=
6\gamma+2$ (resp. $m_{2\gamma+2}(H(P))=6\gamma+5$), and
$\rho(H(P))=2\gamma+1$ (resp. $\rho(H(P))=2\gamma+2$) provided
$\gamma$ odd (resp. $\gamma$ even).
\end{remark}
%
%
%
\section{Weight of semigroups}
\subsection{Bounding the weight.}
Let $H$ be a semigroup of genus $g$. Set $m_j=m_j(H)$ for each $j$ and
$\rho=\rho(H)$ (see (\ref{even-gap}). Due to $m_g=2g$ (see
(\ref{prop-sem})), the weight $w(H)$ of $H$ can be
computed by
\begin{equation}\label{weight}
w(H)=\frac{3g^2+g}{2}-\mathop{\sum}\limits^{g}_{j=1} m_j.
\end{equation}
Consequently the problem of bounding $w(H)$ is equivalent to
the problem of bounding
$$
S(H):= \sum_{j=1}^{g} m_g.
$$
If $\rho=0$, then we have $m_i=2i$ for each $i=1,\ldots,g$. In
particular we have
$w(H)=\binom{g}{2}$. Let $\rho\ge 1$ (or
equivalently $f_1\ge 4$). Then by (\ref{gene}) we have
\begin{equation}\label{weight1}
S(H)=\sum_{f\in \tilde H,\ f\le g} 2f + \sum_{i=1}^{\rho} u_i,
\end{equation}
where
$$
\tilde H:= \{f/2 : f\in H,\ f\ {\rm even}\}.
$$
\begin{lemma'}\label{bounds} With the notation of \S1 we have:
\begin{itemize}
\item[(i)] If $f_1=4$, then $f_i=4i$ for $i=1,\ldots,\rho$.
\item[(ii)] If $f_1\ge 6$, then
$$
f_i\ge 4i+2\ \ {\rm for}\ i=1, \ldots,\rho-2,\ \ f_{\rho-1}\ge
4\rho-4,\ \ f_{\rho}=4\rho.
$$
\item[(iii)] $f_i\le 2\rho+2i$ for each $i$.
\item[(iv)] $2g-4j+1 \le u_j \le 2g-2j+1$, for $j=1,\ldots,\rho$.
\end{itemize}
\end{lemma'}
\begin{proof} By (\ref{feto}), we have
$$
\tilde H= \{\frac{f_1}{2},\ldots, \frac{f_\rho}{2}\}\cup \{4\rho +i:
i\in \mathbb N \}.
$$
Thus $\tilde H$ is a semigroup of genus $\rho$. Then (i) is due to the
fact that $f_1/2=2$ and (ii) follows from (\ref{prop-sem}). Statement (iii)
follows from (\ref{even}).
\smallskip

(iv) The upper bound follows from $u_1\le 2g-1$. To
prove
the lower bound we procced by induction on $i$. The case $i=\rho$ follows
from (\ref{des-odd-1}). Suppose that $u_i\ge 2g-4i+1$ but
$u_{i-1} < 2g-4(i-1)+1$, for $1<i\le \rho$. Then $u_i=2g-4i+1$,
$u_{i-1}= 2g-4i+3$, and
there exists an odd gap $\ell$ of $H$ such that
$\ell>u_{i-1}$. Take the smallest $\ell$ with such a property.
Set $I:=[\ell-u_{i-1}, \ell-u_\rho]\subseteq [2,4\rho-2]$ and let
$t$ be the number of non-gaps of $H$ belonging to $I$. By the election
of $\ell$ we have that $\ell-u_{i-1}<f_1$. Now, since  $\ell-u_j \in I$ for
$j=i-1,\ldots,\rho$ we also have that
$$
\frac{u_{i-1}-u_\rho}{2}+1 \ge t + \rho -(i-1)+1.
$$
Thus $u_\rho \le 2g-2\rho-2i-2t+1$. Now, since
$u_\rho+f_{t+1}>u_{i-1}$ and since by statement (iii) $f_{t+i-1}\le
2\rho + 2t +2i-2$, we have that the odd non-gaps $u_\rho
+f_{t+1}, \ldots, u_\rho +f_{t+i-1}$ belong to $[\ell+2,2g-1]$.
This is a contradiction because $H$ would have $(\rho -i+2)+(i-1) = \rho
+1$ odd non-gaps.
\end{proof}
\begin{lemma'}\label{bo-weight} Let $H$ be a semigroup of genus
$g$. With notation as in \S1, we have
\begin{itemize}
\item[(i)] $w(H)\ge \binom{g-2\rho}{2}$. Equality holds if and
only if $f_1=2\rho+2$ and $u_{\rho}=2g-2\rho+1$.
\item[(ii)] If $g\ge 2\rho$, then $w(H)\le \binom{g-2\rho}{
2}+2\rho^2 $. Equality holds if and only if
 $H=\langle 4, 4\rho,2g-4\rho+1\rangle$.
\item[(iii)] If $g\le 2\rho-1$, then $w(H)\le \binom{g+2\rho}{
2}-4g-6\rho^2+8\rho $.
\end{itemize}
\end{lemma'}
\begin{proof} (i) By (\ref{weight}) we have to show that
$$
S(H) \le
g^2+(2\rho+1)g-2\rho^2-\rho,
$$
and that the equality holds if and only if $f_1=2\rho+2$ and
$u_{\rho}=2g-2\rho+1$. Both the above statements follow from Lemma
\ref{bounds} (i), (iv).
\smallskip

(ii) Here we have to show that
\begin{equation*}
S(H)\ge g^2 + (2\rho+1)g-4\rho^2-\rho,\tag{$\dag$}
\end{equation*}
and that equality holds if and only $H=\langle
4,4\rho+2,2g-4\rho+1\rangle$.

Since $g\ge 2\rho$ by (\ref{feto}) we obtain
\begin{equation}\label{sum1}
S(H)=\sum_{i=1}^{\rho}(f_i(H)+u_i(H))+ g^2+g-4\rho^2-2\rho.
\end{equation}
Thus we obtain $(\dag)$ by means of Lemma \ref{bounds} (ii),
(iii) and (iv). Moreover the equality in $(\dag)$ holds if and only
if $\sum_{i=1}^{\rho}(f_i+u_i)=2\rho g +\rho$. Then the
second part of (ii) also follows from the above mentioned results.
\smallskip

(iii) In this case, due to the proof of Lemma \ref{feto1}, instead
of equation (\ref{sum1}) we have
\begin{equation}\label{sum2}
S(H) = \sum_{i=1}^{2g-3\rho}f_i +
\sum_{i=2g-3\rho+1}^{g-\rho}(2i+2\rho) +
\sum_{i=1}^{\rho}u_i.
\end{equation}
We will see in the next remark that in this case we have $f_1\ge 6$.
Thus by using Lemmas \ref{feto2} and \ref{bounds} (iii), (iv) we
obtain
$$
S(H)\ge 4\rho^2-(2g+7)\rho +g^2 +5g,
$$
from where it follows the proof.
\end{proof}
\begin{remark'}\label{remark3} (i) If $f_1=4$, then $g\ge 2\rho$.
This follows from
the fact the biggest even gap of $H$ is $4(\rho-1)+2$. Moreover,
one can determinate $u_{\rho},\ldots,u_1$ as
follows: let $J\in \mathbb N$ satisfying  the inequalities below
$$
{\rm max}\{1, \frac{3\rho+2-g}{2}\}\le J\le {\rm
min}\{\rho+1,\frac{g-\rho+3}{2}\},
$$
provided $g$ even, otherwise replace $J$ by $\rho-J+2$; then
\begin{eqnarray*}
\{u_{\rho},\ldots,u_1\} & = &
\{ 2g-4\rho+4J-7+4i: i=1,\ldots,\rho-J+1\}\\
&   &\mbox{}\cup\{2g-4J+3+4i: i=1,\ldots,J-1\}.
\end{eqnarray*}
(see \cite[\S3]{Ko}, \cite[Remarks 2.5]{To}). Consequently from
(\ref{sum1}) and (\ref{weight}) we obtain
$$
w(H)=\binom{g-2\rho}{2}+ 2\rho^2+4\rho+6+4J^2-(4\rho+10)J.
$$
In particular we have
$$
\binom{g-2\rho}{2}+\rho^2-\rho\le w(H)\le
\binom{g-2\rho}{2}+2\rho^2.
$$
Let $C$ be an integer such that $0\le 2C\le \rho^2+\rho$. Then
$w(H)=\binom{g-2\rho}{2}+\rho^2-\rho +2C$
if and only if $4+32C$ is a square. The lower bound is attained if
and only if $H=\langle
4,4\rho+2, 2g-2\rho+1,2g-2\rho+3\rangle$.
\smallskip

\noindent (ii) Let $u_{\rho}=3$. Them from (\ref{des-odd-1}) and
Lemma \ref{feto2} we find that $g\in
\{2\rho-1,2\rho,2\rho+1\}$. Moreover, in this case one
can also obtain a explicit formula for $w(H)$ (\cite[Lemma 6]{K1}). Let
$g\equiv r \pmod{3}$, $r=0,1,2$ and let $s$ be an integer such that $0\le
s\le (g-r)/3$. If $r=0,1$ (resp. $r=2$), then
\begin{align*}
w(H) & =\frac{g(g-1)}{3}+3s^2-gs-s\le \frac{g(g-1)}{3} \\
\intertext{resp.}
w(H) & =\frac{g(g-2)}{3}+3s^2-gs+s\le \frac{g(g-2)}{3}.
\end{align*}
If $r=0,1$ (resp. $r=2$), equality occurs if and only if $H=\langle 3,
g+1\rangle$  (resp. $H=\langle 3, g+2,2g+1\rangle$).
\end{remark'}

Let $g\le 2\rho-1$. The way how we bound from below
equation (\ref{sum2}) is far away from being sharp. We do
not know an analogous to the lower bound of Lemma \ref{bounds} (iv) in
this case. However, for certain range of $g$ the bounds in \ref{remark3}
(ii) are the best possible:
\begin{lemma'}\label{opt-weight} Let $H$ be a semigroup of genus
$g\ge 11$, $r\in \{1,2,3,4,5,6\}$ such that $g\equiv r \pmod{6}$. Let
$\rho$ be the number of even gaps of $H$.
If
$$
\rho>\left\{
\begin{array}{ll}
\frac{g-5}{6}  &  {\rm if\ } r=5 \\
\frac{g-r}{6}-1 &  {\rm if\ } r\neq 5,
\end{array}
\right.
$$
then
$$
w(H)\le \left\{
\begin{array}{ll}
\frac{g(g-2)}{3} & {\rm if\ } r= 2,5 \\
\frac{g(g-1)}{3} & {\rm if\ } r=1,3,4,6.
\end{array}
\right.
$$
If $r=2,5$ (resp. $r\not\in\{2,5\}$) equality above holds if and only if
$H=\langle 3,g+2,2g+1\rangle$ (resp. $H=\langle 3,g+1\rangle$).
\end{lemma'}
\begin{proof} We assume $g\equiv 5 \pmod{6}$; the other cases can be
proven in a similar way. By Remark \ref{remark3} (ii) we can assume
$u_1>3$, and then by (\ref{weight}) we have to prove that
\begin{equation*}
S(H) > \frac{7g^2+7g}{6}.\tag{$*$}
\end{equation*}
Now, since $\rho>(g-5)/6$, by Theorem \ref{char4} and Lemma
\ref{feto0} we must have
$$
m_{\frac{g+1}{3}}=m_{\frac{g-1}{2}-\frac{g-5}{6}}\ge g.
$$
(A) Let $S':= \sum_{i} m_i$, $(g+1)/3\le i \le g$:\quad Define
$$
F:= \{ i\in \mathbb N: \frac{g+1}{3}\le i\le g,\ m_i\le 2i+\frac{g-5}{3}\},
$$
and let $f:= {\rm min}(F)$. Then $f\ge (g+4)/3$, $m_f=2f+\frac{g-5}{3},
m_{f-1}=
2f+\frac{g-8}{3}$. Thus for $g\ge i\ge f$, $d_i=1$ and hence by Corollary
\ref{cor-cast}, $2m_i\ge m_{3i-1}=g+3i-1$. In particular, $f\ge (g+7)/3$.
Now we bound $S'$ in three steps:
\smallskip

\noindent Step (i). $(g+1)/3\le i\le f-1$: By definiton of $f$ we have that
$m_i\ge 2i+\frac{g-2}{3}$ and hence
\begin{equation}\label{aux0}
\sum_{i} m_i \ge f^2+\frac{g-5}{3}f -\frac{2g^2-2g-4}{9}.
\end{equation}

\noindent Step (ii). $f\le i\le (6f-g-7)/3$: Here we have that $m_i\ge
m_f+i-f=i+f+\frac{g-5}{3}$. Hence $$
\sum_{i} m_i \ge \frac{5}{2}f^2-\frac{4g+37}{6}f -\frac{g^2-13g-68}{18}.
$$

\noindent Step (iii). $(6f-g-4)/3\le i\le g$: Here we have $m_i+m_{i+1}\ge
g+3i+1$ for $i$ odd, $6f-g-4\le i\le g-2$. Since $m_g=2g$ then we have
$$
\sum_{i} m_i \ge -3f^2+6f+\frac{4g^2+2g-8}{3}.
$$
(B) Let $S'':= \sum_{i} m_i$, $1\le i\le (g-2)/3$:\quad By Theorem
\ref{char2} and Lemma \ref{feto0} we have that $m_i\ge 3i$ for
$i$ odd, $i=3,\ldots, (g-2)/3$. First we notice that for $i$ odd and
$3\le i\le (g-8)/3$ we must have $m_{i+1}\ge 3i+3$. Otherwise we would
have $d_{i+1}=1$ and hence by Corollary \ref{cor-cast} and
(\ref{prop-sem}) we would have $2m_{i+1}\ge m_{3i+2}\ge 6i+5$, a
contradiction.
\begin{claim*} Let $i$ odd and $3\le i\le (g-8)/3$. If
$m_i=3i$ or $m_{i+1}=3i+3$, then $m_1=3$.
\end{claim*}
\begin{proof} {\it (Claim).} It is enough to show that
$d_i=3$ or $d_{i+1}=3$. Suppose that $m_i=3i$. Since
$i$ is odd, $d_i$ is one or three. Suppose $d_i=1$. Then by Corollary
\ref{cor-cast} we have $6i=2m_i\ge m_{3i-1}$ and hence
$6i=2m_i=m_{3i-1}$. Let $\ell \in G(H)$. Then $\ell \ge m_{3i-1}+3$. In
fact if $\ell>m_{3i-1}+3$, by choosing the smallest $\ell$ with such a
property we would have
$3i+2$ gaps in $[1, 6i]$ namely,
$1,2,3,\ell-m_{3i-1},\ldots,\ell-m_1$, a contradiction.
Then it follows that $g\le 3i+1+3=3i+4$ or $g+2\le 3i+4$.
\smallskip

Now suppose that $m_{i+1}=3i+3$; as in the previous proof here we also have
that $d_{i+1}>1$. Suppose that $d_{i+1}=2$. Then $m_1>3$ and hence
$m_i=3i+1$. Since we know that $m_{i+2}\ge 3i+6$, then the even number
$\ell=3i+5$ is a gap of $H$. Then we would fine $2i+2$ even numbers in
$[2,3i+3]$, namely $m_1,\ldots,m_{i+1}$, and $\ell-m_{i+1},\ell-m_1$, a
contradiction. Hence $d_{i+1}=3$ and then $m_1=3$.
\end{proof}
Then, since we assume $u_1>3$, we have $m_i+m_{i+1}\ge 6i+5$ for $i$ odd
$3\le i\le (g-8)/3$, $m_{\frac{g-2}{3}}\ge g-2$, and so
\begin{equation}\label{aux}
\begin{split}
\sum_{i=1}^{(g-2)/3} m_i & \ge \sum_{j=1}^{(g-11)/6} (12j+10)
+m_1+m_2+m_{\frac{g-2}{3}}    \\
             & \ge \frac{g^2+g-78}{6} + m_1+m_2.\\
\end{split}
\end{equation}

Summing up (i), (ii), (iii) and (B) we get
$$
S(H)\ge \frac{3f^2-(2g+11)f}{6}+\frac{22g^2+32g-206}{18}+m_1+m_2.
$$
The function $\Gamma(x):= 3x^2-(2g+11)x$ attains its minimum for
$x=(2g+11)/6<(g+7)/3\le f$. Suppose that $f\ge (g+13)/3$. Then we find
$$
S(H)\ge \frac{7g^2+7g-60}{6}+m_1+m_2.
$$
We claim that $m_1+m_2>11$. Otherwise we would have $m_3=m_1+m_2=10$
which is impossible. From the claim we get $(*)$.

In all the computations below we use the fact that $2g\le (m-1)(n-1)$
whenever $m,n \in H$ with $\gcd(m,n)=1$ (see e.g. Jenkins \cite{J}).

Now suppose that $f=(g+10)/3$. Here we find
$$
S(H)\ge \frac{7g^2+7g-72}{6}+m_1+m_2.
$$
Suppose that $m_1+m_2\le 12$ (otherwise the above computation imply
$(*)$.). If $g>11$,
then $m_4\ge 13$ and so $m_3=m_1+m_2\in \{9,11,12\}$. If $m_1+m_2=9$,
then $g\le 6$; if $m_1+m_2=11$ then $g\le 10$; if $m_1+m_2=12$ then $g\le
11$ or $m_1=4$, $m_2=8$. Let
$s$ denote the first odd non-gap of $H$. Then $2g\le 3(s-1)$ and so
$s>(2g+2)/3$. In the interval $[4,(2g+2)/3]$ does not exist $h\in H$ such
that $h\equiv 2 \pmod{4}$: In fact if such a $h$ exists then we would
have $4\rho+2\le (2g-4)/3$ or $ \rho\le (g-5)/6$. Consequently
$m_3=12,\ldots,m_{(g+1)/6}=(2g+2)/3$. Thus we can improve the
computation in (\ref{aux}) by summing it up $\sum_{i=1}^{j}(4i+1)$, where
$j=(g-5)/12$ or $j=(g-11)/12$. Then we get $(*)$. If $g=11$, the first
seven non-gaps are $\{4,8,10,12,14,15,16\}$ or $\{5,7,10,12,14,15,16\}$.
In both cases the computation in (\ref{aux0}) increases at least by one,
and so we obtain $(*)$.

Finally let $f=(g+7)/3$. Here we find
$$
S(H)\ge \frac{7g^2+7g-78}{6}+m_1+m_2,
$$
and we have to analize the cases $m_1+m_2=\{9, 11, 12, 13\}$. This can be
done as in the previous case. This finish the proof of Lemma
\ref{opt-weight}.
\end{proof}
\subsection{The equivalence $(P_1)\Leftrightarrow (P_4)$.}
We are going to characterize $\gamma$-hyperelliptic semigroups by means
of weights of semigroups. We begin with the following result, which
has been proved by Garcia for $\gamma\in\{1,2\}$ \cite[Lemmas 8 and 10]{G}.
\begin{theorem'}\label{char-weight} Let $\gamma\in \mathbb N$ and $H$ a
semigroup whose genus $g$ satisfies
(\ref{bound-g}). Then the following statements are equivalent:
\begin{itemize}
\item[(i)] $H$ is $t$-hyperelliptic for some $t\in \{0,\ldots,\gamma\}$.
\item[(ii)] $w(H)\ge \binom{g-2\gamma}{2}$.
\end{itemize}
\end{theorem'}
\begin{theorem'}\label{char-weight1} Let $\gamma$, $H$ and $g$ be
as in
Theorem \ref{char-weight}. The following statements are equivalent:
\begin{itemize}
\item[(i)] $H$ is $\gamma$-hyperelliptic.
\item[(ii)] $\binom{g-2\gamma}{2}\le w(H)\le \binom{g-2\gamma}{2}+2\gamma^2$.
\item[(iii)] $\binom{g-2\gamma}{2}\le w(H)<\binom{g-2\gamma+2}{2}$.
\end{itemize}
\end{theorem'}
\begin{proof} {\it (Theorem \ref{char-weight}).} (i) $\Rightarrow$ (ii):
By Lemma \ref{feto0} and
Lemma \ref{bo-weight} (i) we have $w(H)\ge \binom{g-2t}{2}$. This implies
(ii).

(ii) $\Rightarrow$ (i): Suppose that $H$ is not $t$-hyperelliptic for any
$t\in \{0,\ldots\gamma\}$. We are going to prove that
$w(H)<\binom{g-2\gamma}{2}$, which
by (\ref{weight}) is equivalent to prove that:
\medskip

$(*)$\hfill $\sum_{i=1}^{g} m_i >
g^2+(2\gamma+1)g-2\gamma^2-\gamma.$\hfill
\medskip

\noindent We notice that by Lemma \ref{feto0} we
must have $\rho\ge \gamma+1$.
\smallskip

\noindent Case 1: $g$ satisfies the hypothesis of Lemma
\ref{opt-weight}.\quad From that lemma we have $ S(H)\ge (7g^2+5g)/3$
and then we get $(*)$ provided
$$
g>\bar\gamma:= 12\gamma+1+\sqrt{96\gamma^2+1}.
$$
We notice that
$\gamma^2+4\gamma+3\ge \bar\gamma$ if $\gamma\ge 7$. For $\gamma=1,4,6$ we
need respectively $g>11$, $g>44$ and $g>65$. By noticing that 11, 44 and
65 are
congruent to 2 modulus 3, we can use $g=11$, $g=44$ and $g=65$ because in
these cases
$S(H)\ge (7g^2+7g)/3$. For the other values of $\gamma$ we obtain the
bounds of (\ref{bound-g}).
\smallskip

\noindent Case 2: $g$ does not satisfy the hypothesis of Lemma
\ref{opt-weight}.\quad Here we have $g\ge 6\rho+5$. From (\ref{sum1})
and Lemma \ref{bounds} we have $S(H)\ge
g^2+(2\rho+1)g-4\rho^2-\rho$. The function
$\Gamma(\rho):= (2g-1)\rho-4\rho^2$ satisfies
$$
\Gamma(\rho)\ge \Gamma(\gamma+1)=(2\gamma+2)g-4\gamma^2-9\gamma-5,
$$
for $\gamma+1 \le \rho\le [(2g-1)/4]-\gamma-1$. Thus we obtain
condition $(*)$ provided $g\ge \gamma^2 + 4\gamma+3$.
\end{proof}
\begin{remark'}\label{oliv} Let $H$ be a semigroup of genus $g$, $r$ the
number
defined in Lemma \ref{opt-weight}. Put $c:= (g-5)/6$ if $r=5$, and $c:=
(g-r)/6-1$ otherwise.

{}From the proof of Case 2 of the above result we see that $S(H)\ge
g^2+3g-5$ whenever $1\le \rho(H)\le (g-3)/2$. Hence this result and Lemma
\ref{opt-weight} imply
$$
w(H)\le\left\{
\begin{array}{ll}
(g^2-5g+10/2 & {\rm if\ } \rho(H)\le c\\
{\rm min}\{(g^2-5g+10)/2, (g-1)g/3\}  &  {\rm if\ }
c<\rho(H)\le (g-3)/2\\
(g-1)g/3    &  {\rm if\ } \rho(H)>(g-3)/2.
\end{array}
\right.
$$
\end{remark'}
\begin{proof} {\it (Theorem \ref{char-weight1}).} (i)
$\Rightarrow$ (ii) follows from Lemma \ref{bo-weight}. (ii) $\Rightarrow$
(iii) follows from the hypothesis on $g$.

(iii) $\Rightarrow$ (i): By Theorem \ref{char-weight} we have that
$H(P)$ is $t$-hyperelliptic for some $t\in \{0,\ldots,\gamma\}$. Then by
Lemma \ref{bo-weight} and hypothesis we have
$$
\binom{g-2\gamma+2}{2}>w(H)\ge \binom{g-2t}{2},
$$
from where it follows that $t=\gamma$.
\end{proof}
\begin{remark'}\label{sharp-weight} The hypothesis on $g$ in the
above two theorems is sharp:
\smallskip

(i) Let $\gamma\ge 7$ and considerer $H:=\langle 4, 4(\gamma+1),
2g-4(\gamma+1)+1\rangle$ where $g$ is an integer satisfying ${\rm
max}\{4\gamma+4, \frac{\gamma^2+6\gamma-3}{2}\}<g\le \gamma^2+4\gamma+2$.
Then $H$ has genus $g$ and $\rho(H)=\gamma+1$. In particular $H$
is not $\gamma$-hyperelliptic. By Lemma \ref{bo-weight} (ii) we have
$w(H)=\binom{g-2(\gamma+1)}{2}+2(\gamma+1)^2$. Now
it is easy to check that $w(H)$ satisfies Theorem \ref{char-weight} (ii)
and Theorem \ref{char-weight1} (iii).
\smallskip

(ii) Let  $\gamma\le 6$ and consider $H=\langle 3,g+1\rangle$,
where $g=10, 22, 33, 43, 55, 64$ if $\gamma=1,2,3,4,5,6$ respectively.
$H$ has genus $g$ and it can be easily checked that $H$ is not
$\gamma$-hyperelliptic by means of inequality (\ref{des-odd-1}) and Lemma
\ref{feto1}. Moreover
$w(H)=g(g-1)/3$ (see Remark \ref{remark3} (ii)). Now it is easy to check
that $w(H)$ satisfies Theorem \ref{char-weight} (ii) and Theorem
\ref{char-weight} (iii).
\smallskip

(iii) The semigroups considered in (i) and (ii) are also Weierstrass
semigroups (see Komeda \cite{Ko}, Maclachlan \cite[Thm. 4]{Mac}).
\end{remark'}
\subsection{Weierstrass weights}
In this section we apply Theorem \ref{char-weight1} in order to
characterize double coverings of curves by means of Weierstrass weights.
Specifically we strengthen \cite[Theorem B]{To} and hence all its
corollaries. The basic references for the discussion below are Farkas-Kra
\cite[III.5]{F-K}
and St\"ohr-Voloch \cite[\S1]{S-V}.

Let $X$ be a non-singular,
irreducible and projective algebraic curve of
genus $g$ over an algebraically closed field $k$ of characteristic $p$.
Let $\pi:X\to \mathbb P^{g-1}$ be the morphism induced by the canonical
linear system on $X$. To any $P\in X$ we can associate the sequence
$j_i(P)$ ($i=0,\ldots,g-1$) of intersection multiplicities at $\pi(P)$ of
$\pi(X)$ with hyperplanes of $\mathbb P^{g-1}$. This sequence is the same
for all but finitely many points (the so called Weierstrass points of
$X$). These points are supported by a divisor $\mathcal{W}$ in such a way
that the Weierstrass weight at $P$, $v_P(\mathcal{W})$, satisfies
$$
v_P({\mathcal{W}})\ge w(P):= \sum_{i=1}^{g-1}(j_i(P)-\epsilon_i),
$$
where $0=\epsilon_0<\ldots<\epsilon_{g-1}$ is the sequence at a generic
point. One has $j_i(P)\ge \epsilon_i$ for each $i$, and from the
Riemann-Roch theorem follows that $G(P):=\{j_i(P)+1:i=0,\ldots,g-1\}$ is
the set of gaps of a semigroup $H(P)$ of genus $g$ (the so called
Weierstrass
semigroup at $P$). $X$ is called {\it classical} if $\epsilon_i=i$ for
each $i$ (e.g. if $p=0$ or $p>2g-2$). In this case we have
$v_P({\mathcal{W}})=w_P(R)$
for each $P$, and the number $w(P)$ is just the weight of the
semigroup $H(P)$ defined in
\S0. The following result strengthen \cite[Thm.B]{To}. The proof follows
from \cite[Thm.A]{To}, \cite[Thm.A]{To1}, and Theorem \ref{char-weight1}.
\begin{theorem'} Let $X$ be a classical curve, and assume that $g$
satisfies (\ref{bound-g}). Then the following statements are equivalent:
\begin{itemize}
\item[(i)] $X$ is a double covering of a curve of genus $\gamma$.
\item[(ii)] There exists $P\in X$ such that
$$
\binom{g-2\gamma}{2}\le w(P)\le \binom{g-2\gamma}{2}+2\gamma^2.
$$
\item[(iii)] There exists $P\in X$ such that
$$
\binom{g-2\gamma}{2}\le w(P)< \binom{g-2\gamma+2}{2}.
$$
\end{itemize}
\end{theorem'}
Remark \ref{sharp-weight} says that the bound for $g$ above is the best
possible. Further applications of \S3.1 and \S3.2 will be published
elsewhere \cite{To2}.

\end{document}